\documentclass[a4paper, 11pt]{amsart}

\usepackage{amssymb,amsmath,latexsym,amsthm}
\usepackage[latin1]{inputenc}
\usepackage[T1]{fontenc}
\usepackage[francais,english]{babel}
\usepackage{amssymb}
\usepackage{amsmath}
\usepackage{amsthm}
\usepackage{amscd}
\usepackage{amsfonts}
\usepackage{stmaryrd}
\usepackage{pb-diagram}
\usepackage{epic,eepic,epsfig}
\usepackage{a4wide}
\usepackage{nextpage}
\usepackage{fancyhdr}
\newtheorem{prop}{Proposition}[section]
\newtheorem{cor}[prop]{Corollaire}

\newtheorem{lem}[prop]{Lemme}
\newtheorem{thm}[prop]{Th{\'e}or{\`e}me}
\newtheorem{conj}[prop]{Conjecture}

\newtheorem{fait}[prop]{Fait}
\newtheorem{theorem}{Th{\'e}or{\`e}me}
\newtheorem{corol}{Corollaire}

\newtheorem{Conjecture}{Conjecture}

\newtheorem*{remarque}{Remarque}
\newtheorem*{exemple}{Exemple}
\newtheorem*{applications}{Applications}

\newcommand{\Jac} {\mathop{\mathrm{Jac}}}
\newcommand{\Ima} {\mathop{\mathrm{Im}}}

\newcommand{\disc} {\mathop{\mathrm{disc}}}

\newcommand{\Ree} {\mathop{\mathrm{Re}}}

\newcommand{\Divi} {\mathop{\mathrm{Div}}}

\newcommand{\Nk} {\mathop{\mathrm{N}_{k/\mathbb{Q}}}}

\newcommand{\dime} {\mathop{\mathrm{dim}}}

\newcommand{\NL} {\mathop{\mathrm{N}_{L/k}}}
\newcommand{\Pico} {\mathop{\mathrm{Pic^{0}}}}
\newcommand{\NkN} {\mathop{\mathrm{N}_{k_{N}/k}}}
\newcommand{\hF} {\mathop{h_{\mathrm{F}}}}

\newcommand{\hst} {\mathop{h_{\mathrm{st}}}}

\newcommand{\Gal} {\mathop{\mathrm{Gal}}}
\newcommand{\Aut} {\mathop{\mathrm{Aut}_{\mathbb{F}_{l}}}}

\begin{document}

\def\refname{\centerline{Bibliographie}}

\title[Conjecture de Lang]{Remarques sur une conjecture de Lang}


\author[Fabien {\sc Pazuki}]{{\sc Fabien} Pazuki}
\address{Fabien {\sc Pazuki}\\
IMJ Universit\'e Paris 7\\
2, place de Jussieu\\
75 251 Paris Cedex 05, France}
\email{pazuki@math.jussieufr}
\urladdr{http~://www.math.jussieu.fr/~pazuki}

\maketitle
\begin{centering}
{\small{\textsc{R\'esum\'e} : Le but de cet article est d'étudier une conjecture de Lang énoncée sur les courbes elliptiques dans un livre de Serge Lang, puis généralisée aux variétés abéliennes de dimension supérieure dans un article de Joseph Silverman. On donne un résultat asymptotique sur la hauteur des points de Heegner sur $J_{0}(N)$, lequel permet de déduire que la conjecture est optimale dans sa formulation.}}
\end{centering}

\begin{abstract}
The aim of this paper is to study a conjecture predicting a lower bound on the canonical height on abelian varieties, formulated by S. Lang and generalized by J. H. Silverman. We give here an asymptotic result on the height of Heegner points on the modular jacobian $J_{0}(N)$, and we derive non-trivial remarks about the conjecture.
\end{abstract}

\bigskip

\section{La conjecture de Lang et Silverman}

S. Lang a conjectur{\'e} dans \cite{Lan} p. 92 une minoration de la
hauteur de N{\'e}ron-Tate d'une courbe elliptique, qu'on rappelle ici~:

\begin{Conjecture}(Lang)
Pour tout corps de nombres $k$, il existe une constante positive $c(k)$
telle que pour toute courbe elliptique $E$ d{\'e}finie sur $k$ et tout
point $P$ d'ordre infini de $E(k)$ on ait~:
\[
 \widehat{h}(P) \geq c(k)\, \max\Big\{\log \Nk(\Delta_{E}),h(j_{E})\Big\},
\]
o{\`u} $\widehat{h}(.)$ est la hauteur de N{\'e}ron-Tate sur $E$,
$\Nk(\Delta_{E})$ la norme du disciminant minimal de la
courbe $E$ et $h(j_{E})$ la hauteur de Weil logarithmique et absolue
de l'invariant modulaire $j_{E}$ de la courbe $E$.

\end{Conjecture}

\begin{remarque}
 Dans cette conjecture il est équivalent de chercher une minoration du
type $\widehat{h}(P)\geq c(k)\, \hF(E/k)$ où $\hF(E/k)$ est la hauteur de Faltings (relative)
de la courbe elliptique $E$. Dans la formulation de la question qui figure dans \cite{Lan}, S. Lang ne faisait intervenir que le logarithme du discriminant.
\end{remarque}

Cette conjecture de Lang a {\'e}t{\'e} partiellement d{\'e}montr{\'e}e par
M. Hindry et J. Silverman qui obtiennent dans \cite{HiSi3}, corollaire 4.2 (ii) de leur théorème 4.1 (p. 430 et 431), le résultat suivant~:

\begin{theorem}(Hindry, Silverman) \label{elliptique}
Soit $k$ un corps de nombres de degré $d$. Soit $E/k$ une courbe elliptique de
disciminant minimal $\Delta_{E}$ et de conducteur $F_{E}$. On note $\sigma_{E}$ le
quotient de Szpiro défini par $\sigma_{E} = \log \Nk(\Delta_{E})/\log \Nk(F_{E})$. Alors
pour tout point $P\in{E(k)}$ d'ordre infini on a la minoration~:
\[
\widehat{h}(P)\geq (20\sigma_{E})^{\displaystyle{-8d}}10^{\displaystyle{-4\sigma_{E}}}\frac{1}{12}\max\Big\{\log \Nk(\Delta_{E}),h(j_{E})
\Big\}.
\]
\end{theorem}

 Ceci permet de conclure pour toute famille de courbes elliptiques pour lesquelles le
quotient de Szpiro est borné uniformément. Une conjecture de Szpiro affirme que
c'est en fait le cas de toutes les courbes elliptiques sur $k$ et entraîne donc la
conjecture de Lang ci-dessus. La preuve de ce théorème repose sur l'existence d'une décomposition de la hauteur de Néron-Tate en somme de hauteurs locales bien normalisées.

J. Silverman avait démontré auparavant plusieurs cas particuliers de cette conjecture dans \cite{Sil4} et \cite{Sil3}. Par la suite S. David a publié une preuve de transcendance \cite{Dav3} offrant une constante $c(d,\sigma_{E})$ polynomiale inverse en $d$ et $\sigma_{E}$. On peut citer aussi l'article de M. Krir
\cite{Krir} qui explicite sur $k=\mathbb{Q}$ d'une manière un peu différente ce résultat de minoration pour des familles de courbes elliptiques particulières. Plus récemment, une nouvelle constante polynomiale inverse a été obtenue par C. Petsche \cite{Pet} par la technique de décomposition locale.

\vspace{0.3cm}

La conjecture sur les courbes elliptiques a ensuite été g{\'e}n{\'e}ralis{\'e}e aux vari{\'e}t{\'e}s ab{\'e}liennes de dimension sup{\'e}rieure par J. Silverman dans \cite{Sil3} p. 396~:

\begin{Conjecture}(Lang, Silverman)\label{Lang Silverman}
Soit $g \geq 1$. Pour tout corps de nombres $k$, il existe une
constante positive $c(k,g)$ telle que pour toute vari{\'e}t{\'e} ab{\'e}lienne
$A/k$ de dimension $g$, pour tout diviseur ample et symétrique $\mathcal{D}\in{\Divi(A)}$ et
tout point $P\in{A(k)}$ tel que $\mathbb{Z}\!\cdot\! P=\{mP|m\in{\mathbb{Z}}\}$ soit Zariski-dense on ait~:
\[
 \widehat{h}_{A,\mathcal{D}}(P) \geq c(k,g)\, \max\Big\{1,\hF(A/k)\Big\},
\]
o{\`u} $\widehat{h}_{A,\mathcal{D}}(.)$ est la hauteur de N{\'e}ron-Tate sur $A$ associ{\'e}e au
diviseur $\mathcal{D}$ et $\hF(A/k)$ est la hauteur de Faltings (relative) de la vari{\'e}t{\'e}
ab{\'e}lienne $A$. 

\end{Conjecture}

\begin{remarque}
Il y a plusieurs notions de hauteur d'une variété abélienne. L'énoncé de cette conjecture est plus fin avec la hauteur de Faltings (relative) comme minorant qu'avec la hauteur de Faltings \emph{stable} notée $h_{\mathrm{st}}$. Rappelons de plus que la hauteur de Faltings stable est comparable à une hauteur modulaire, comme par exemple la hauteur thêta d'une variété abélienne.
\end{remarque}

\begin{remarque}
On peut se demander s'il est possible de conjecturer encore mieux en imposant $c(k,g)=c_{0}$ une constante absolue. On va voir dans cet article que c'est impossible.
\end{remarque}

S. David a proposé une preuve partielle de cette conjecture
g{\'e}n{\'e}ralis{\'e}e, preuve basée sur un raisonnement de type transcendance
(voir \cite{Dav2})~: il donne une borne inf{\'e}rieure pouvant tendre vers l'infini avec la hauteur (thêta) de la vari{\'e}t{\'e}. Plus précisément il obtient le théorème~:

\begin{theorem}(David)
Soient $g\geq 1$ un entier, $k$ un corps de nombres, $v$ une place archimédienne, $(A,\mathcal{D})/k$ une variété abélienne principalement polarisée de dimension $g$ et $\tau_{v}$ une matrice telle que $A(\bar{k}_{v})\cong \mathbb{C}^{g}/\mathbb{Z}^{g}+\tau_{v}\mathbb{Z}^{g}$. On note $||\Ima\tau_{v}||=\max_{i,j} |\Ima\tau_{v,ij}|$. Posons~: $\rho(A)=h_{\mathrm{st}}(A)/\parallel\Ima\tau_{v}\parallel.$

Alors il existe une constante $c_{1}(k,g)>0$ telle que, tout point $P\in{A(k)}$ vérifiant que $\mathbb{Z}.P$ est Zariski-dense, on a~:
\[
 \widehat{h}_{A,\mathcal{D}}(P)\geq c_{1}(k,g)\rho(A)^{-4g-2}\Big(\log\rho(A)\Big)^{-4g-1}\,h_{\mathrm{st}}(A).
\]
\end{theorem}

Cet énoncé implique donc l'inégalité cherchée pour les familles de variétés abéliennes vérifiant $\rho(A,k)$ borné. D. Masser utilise d'ailleurs ces
r{\'e}sultats dans \cite{Masser} pour exhiber une famille de variétés abéliennes simples
avec $\rho$ borné, famille vérifiant donc la conjecture de Lang et Silverman. On trouvera des énoncés plus récents traitant notamment de familles en dimension 2 dans le chapitre 2 de \cite{Paz}.

\vspace{0.3cm}

\begin{applications}
Un résultat de minoration uniforme en la variété du type de l'énoncé de Lang et Silverman aurait des conséquences intéressantes pour plusieurs problèmes concernant les variétés algébriques. On se limitera ici à deux problèmes applicatifs, en direction desquels on trouvera dans la suite des énoncés partiels. Tout d'abord les techniques de preuve des résultats partiels en direction de l'inégalité de Lang et Silverman passent généralement par un raisonnement du type~: \og parmi les $N$ points distincts $P_{1}$,...,$P_{N}$, il en existe un qui vérifie $\widehat{h}(P_{i})>\alpha$ \fg. Si $\alpha$ est strictement positif, on déduit donc qu'il ne peut y avoir plus de $N$ points de hauteur nulle, ce qui procure une borne uniforme sur la torsion des variétés abéliennes considérées pour peu que $N$ soit uniforme. Le deuxième problème lié à ces minorations est l'obtention de bornes uniformes sur le nombre de points rationnels d'une courbe algébrique de genre $g\geq 2$, en passant par l'étude de la variété jacobienne.
\end{applications}

Nous réunissons ici des remarques concernant la conjecture \ref{Lang Silverman}. On montre en particulier qu'il est impossible de proposer une conjecture plus générale dans laquelle la constante de comparaison des hauteurs ne dépend pas du corps ou ne dépend pas de la dimension de la variété. On traite
en détail le cas des jacobiennes de courbes modulaires
$J_{0}(N)$. Plus exactement on produit un équivalent de la hauteur de
Néron-Tate d'un point de Heegner lorsque le niveau $N$ est grand,
généralisant une démarche déjà présente dans \cite{MUll}. Pour $k$ un
corps de nombres dont l'anneau des entiers est noté $\mathcal{O}_{k}$,
on note $h_{k}$ son nombre de classes, $u_{k}$ la moitié du cardinal 
de ses unités et $\mathbb{N}_{k}$ l'ensemble des entiers $N$ tels qu'il existe un point de Heegner associé à $\mathcal{O}_{k}$ sur $X_{0}(N)$. Cet ensemble peut aussi être défini par des congruences. Soulignons de plus qu'on imposera toujours aux entiers considérés dans $\mathbb{N}_{k}$ d'être premiers à $6$ et sans facteur carré.
Le résultat est le suivant~:

\begin{theorem}\label{equiv}
Soit $k$ un corps quadratique dont le discriminant $D$ v\'erifie
les conditions $D<0$ et $D\equiv 1\pmod{4}$. Soit $N\in \mathbb{N}_{k}$, soit $x_{D}\in X_{0}(N)$ un point de Heegner associ\'e \`a $k$ et posons~: $c_{D}=(x_{D})-(\infty)$. Alors on a~:
\[
 \widehat{h}_{J_{0}(N)}(c_{D})\sim h_{k}u_{k}\log (N) ,
\]
lorsque $N\in \mathbb{N}_{k}$ tend vers l'infini.
\end{theorem}

Notons $g(N)$ la dimension de $J_{0}(N)$. L'utilisation de l'équivalent (obtenu dans \cite{JoKr} grâce à des calculs de géométrie hyperbolique complexe) $\hst(J_{0}(N))\sim g(N)\log (N)/3$ de la hauteur de Faltings stable de $J_{0}(N)$ lorsque $N$ est grand et sans facteur carré permet, par comparaison des asymptotiques, de conclure au fait suivant~:

\vspace{0.3cm}

\begin{corol}\label{asympt}
Soit $k$ un corps quadratique dont le discriminant $D$ v\'erifie
les conditions $D<0$ et $D\equiv 1\pmod{4}$. Soit $N\in \mathbb{N}_{k}$, soit $x_{D}\in X_{0}(N)$ un point de Heegner associ\'e \`a $k$ et posons~: $c_{D}=(x_{D})-(\infty)$. Notons $g(N)$ le genre de $X_{0}(N)$. Alors on a~:
\[
\widehat{h}_{J_{0}(N)}(c_{D})\sim \frac{3h_{k}u_{k}}{g(N)}\hst(J_{0}(N)), 
\]
lorsque $N\in \mathbb{N}_{k}$ tend vers l'infini.

\end{corol}

\vspace{0.2cm}

Merci à l'arbitre de la publication qui par ses remarques précises a permis d'améliorer le texte en plusieurs endroits.

\vspace{0.2cm}

\section{Points de Heegner et courbes modulaires}

On s'intéresse dans cette partie aux jacobiennes de courbes modulaires et aux points particuliers que sont les points de Heegner sur ces jacobiennes.

\vspace{0.3cm}

Soit $k$ un corps quadratique imaginaire dont le discriminant $D<0$ est tel que $D\equiv 1 \pmod{4}$ ($D$ sans facteur carré). On s'int\'eresse dans un premier temps \`a l'ensemble $\mathbb{N}_{k}$ des entiers $N$
tels qu'il existe un point de Heegner associ\'e \`a $k$ sur la courbe modulaire $X_{0}(N)$.
On estime ensuite, pour de tels $N$, la hauteur du point de Heegner sur la jacobienne
$J_{0}(N)=\Jac(X_{0}(N))$. On montre en étudiant les différents termes présents l'asymptotique du théorème \ref{equiv}, en notant $h_{k}$ le nombre de classes
associ\'ees \`a $k$ et $u_{k}$ la moiti\'e du nombre de ses unit\'es.

On en d\'eduit ensuite des cons\'equences sur la conjecture de
Lang et Silverman et sur la torsion des jacobiennes de courbes modulaires.

\subsection{Cadre général}

\subsubsection{La courbe modulaire $X_{0}(N)$}
Soit $\mathcal{H}=\{z\in\mathbb{C}\,|\,\Ima(z)>0\}$ le demi-plan de
Poincar\'e et $\mathcal{H}^{*}=\mathcal{H}\cup \mathbb{Q}\cup
\{\infty\}$. Le groupe
\[
\Gamma_{0}(N)=\left\{\left(\begin{array}{cc}
   a & b \\
   c & d \
 \end{array}\right)\in SL_{2}(\mathbb{Z})\,\Big|\,c\equiv 0 \pmod{N}
\right\}
\]
agit sur $\mathcal{H}^{*}$ et
$\mathcal{H}^{\ast}/\Gamma_{0}(N)$ est une surface de Riemann
compacte ; c'est la compactifi{\'e}e de $Y_{0}(N):=
\mathcal{H}/\Gamma_{0}(N)$ laquelle param{\`e}tre les paires
$(E_{\tau},G_{\tau})$ o{\`u} $E_{\tau}=\mathbb{C}/(\mathbb{Z}+\tau
\mathbb{Z})$ est une courbe elliptique sur $\mathbb{C}$ et
$G_{\tau}$ est un sous-groupe de $E_{\tau}(\mathbb{C})$ cyclique
d'ordre $N$. Sur un corps de caract\'eristique $0$ un point de
$Y_{0}(N)$ correspond \`a une paire $(E,E')$ de courbes
elliptiques munies d'une isog\'enie $\phi:E\rightarrow E'$ cyclique de degré $N$.

\vspace{0.3cm}

 Cela permet d'identifier $\mathcal{H}^{\ast}/\Gamma_{0}(N)$ {\`a} la courbe
modulaire $X_{0}(N)$ d{\'e}finie sur $\mathbb{Q}$.

\vspace{0.3cm}

\subsubsection{Points de Heegner}

Soit $k$ un corps quadratique imaginaire. Soit $N$ un entier premier au discriminant de $k$. Le point $x=(E\rightarrow
E')\in X_{0}(N)$ est appel\'e point de Heegner associé à $k$ lorsque les deux
courbes elliptiques $E$ et $E'$ sont \`a multiplication complexe par $\mathcal{O}_{k}$.

On peut d\'ecrire les points de Heegner sur $\mathbb{C}$ (on
pourra se r\'ef\'erer par exemple \`a \cite{GZ} p. 235 et
\cite{Gro})~:

\begin{align*} 
\Big\{(\mathcal{A},\mathbf{n}),\mathcal{A}\in Cl_{k},\mathbf{n}\subset \mathcal{O}_{k}, \mathcal{O}_{k}/\mathbf{n}\cong \mathbb{Z}/N\mathbb{Z} \Big\} &\leftrightarrow \Big\{x_{\mathrm{Heegner}} \in X_{0}(N)(\mathbb{C})\Big\}\\
([\mathbf{a}],\mathbf{n}) & \to (\mathbb{C}/\mathbf{a}\rightarrow \mathbb{C}/\mathbf{a}\mathbf{n}^{-1}).
\end{align*}

La condition d'existence d'un point de Heegner est donc
l'existence d'un id\'eal $\mathbf{n}\subset \mathcal{O}_{k}$ tel
que $\mathcal{O}_{k}/\mathbf{n}\cong \mathbb{Z}/N\mathbb{Z}$. Un
tel id\'eal existe si et seulement s'il existe
$\beta\in\mathbb{Z}/2N\mathbb{Z}$ v\'erifiant $\beta^{2}\equiv D
\pmod{4N}$. On a alors
$\mathbf{n}=\mathbb{Z}N+\mathbb{Z}\frac{\beta +\sqrt{D}}{2}$.

\vspace{0.2cm}

Un point $\tau \in \mathcal{H}$ correspondant \`a une courbe \`a
multiplication complexe par $\mathcal{O}_{k}$ est racine d'une
\'equation quadratique de la forme $A\tau^{2}+B\tau+C=0$ avec $A$,
$B$ et $C$ entiers et de discriminant $B^{2}-4AC=D=\disc(k)$. Or
$N\tau$ doit avoir la m\^eme propri\'et\'e. Ceci implique que
d'une part $N|A$, d'autre part $B^{2}\equiv D \pmod{4N}$. \`A toute
classe de formes quadratiques de discriminant $\disc(k)$ correspond
un point de Heegner diff\'erent associ\'e au corps $k$.

\subsubsection{L'ensemble $\mathbb{N}_{k}$}
On consid\`ere un discriminant de la forme $D=-d_{1}...d_{r}$
avec les $d_{i}$ des entiers naturels premiers impairs deux \`a deux
distincts. Le corps $k$ donnera lieu \`a des points de Heegner sur une courbe
$X_{0}(N)$ pour les entiers $N \in \mathbb{N}_{k}$ avec~:
\[
\mathbb{N}_{k}:=\left\{N\in \mathbb{N}\;\Big|\;
\begin{array}{l}
(1)\; \exists\, d\in{\mathbb{N}},\,d\wedge 4N=1,\, D \equiv d^{2} \pmod{4N}\\
(2)\; N\, \textrm{est sans facteur carr\'e et premier \`a}\; 6
\end{array}
\right\}.
\]
L'hypothèse $(2)$ est essentiellement technique et figure ici pour pouvoir utiliser des calculs asymptotiques plus faciles à mener sous ces hypothèses.

\begin{prop}\label{NK}
L'ensemble $\mathbb{N}_{k}$ est de cardinal infini.
\end{prop}

\begin{proof} Il suffit de traiter le cas $N=p$ premier supérieur ou égal à $5$.
Prenons donc $p\geq 5$ un nombre premier diff\'erent des $d_{i}$. Le discriminant $D$
doit \^etre un carr\'e inversible modulo $4p$, ce qui est \'equivalent aux
conditions $D\equiv 1 \pmod{4}$ et $\left(D/p\right)=1$. 
On obtient pour la seconde~:
\[
\left(\frac{D}{p}\right)=\left(\frac{-1}{p}\right)\prod_{i=1}^{r}\left(\frac{d_{i}}{p}\right)=1.
\]
Une telle \'equation en $p$ admet comme ensemble type de solutions
un nombre non nul et fini de progressions arithm\'etiques (on pourra consulter \cite{IR} p. 55), ceci \'etant une application directe
r\'ep\'et\'ee de la loi de r\'eciprocit\'e quadratique. Le
cardinal de $\mathbb{N}_{k}$ est donc infini.
\end{proof}

\begin{exemple}
On peut traiter un exemple simple pour illustrer
ce propos. Si $D=-3$ on a tout d'abord $D\equiv 1 \pmod{4}$ et de
plus pour $p\geq 5$ premier~:

\[
\left(\frac{D}{p}\right)=\left(\frac{-1}{p}\right)\left(\frac{3}{p}\right)=(-1)^{\displaystyle{\frac{p-1}{2}}}\left(\frac{p}{3}\right)(-1)^{\displaystyle{\frac{p-1}{2}\frac{3-1}{2}}}=\left(\frac{p}{3}\right).
\]

On a donc comme solution la progression arithm\'etique $\{p\equiv
1 \pmod{3} \}\subset \mathbb{N}_{k}$. Ceci permet d'ailleurs de
donner une minoration de la densit\'e $d(\mathbb{N}_{k})$ des premiers de
$\mathbb{N}_{k}$ (au sens de Dirichlet) gr\^ace \`a la forme forte
du th\'eor\`eme de progression arithm\'etique de Dirichlet~:
$d(\mathbb{N}_{k})\geq \frac{1}{\varphi(3)}=\frac{1}{2}$ avec
$\varphi$ l'indicateur d'Euler.
\end{exemple}

\subsection{La formule de Gross-Zagier}

\subsubsection{Accouplement global}
On va reprendre ici l'expression de l'accouplement global des points de
Heegner sur $J_{0}(N)\!\times\! J_{0}(N)$ obtenu par Gross et Zagier
dans \cite{GZ} page 307 et valable pour $(m,N)=1$. On rappelle que
les calculs men\'es dans l'article \cite{GZ} sont faits place par
place, mais qu'\textit{a priori} le symbole $<c,c>_{v}$ n'est pas bien
d\'efini. On peut cependant calculer $<c,d>_{v}$ avec $c\neq d$ et
utiliser le fait que globalement $<c,d>=<c,c>$ car $c-d$ est de
torsion (c'est ici qu'on applique le th\'eor\`eme de
Manin-Drinfeld). C'est la d\'emarche qu'adoptent B. Gross et
D. Zagier dans leur article.

\vspace{0.3cm}

On rappelle ici le cadre dans lequel on se place~:

\vspace{0.3cm}

\begin{itemize}

\item[$\bullet$]{$x_{D}$ \'etant une coordonn\'ee d'un point de Heegner associ\'e
au corps quadratique imaginaire $k=\mathbb{Q}[\sqrt{D}]$, on consid\`ere le point
$c_{D}=(x_{D})-(\infty)\in J_{0}(N)(H)$ et le point $d_{D}=(x_{D})-(0) \in J_{0}(N)(H)$
avec $H$ le corps de classe de Hilbert associ\'e \`a $k$.}

\vspace{0.3cm}

\item[$\bullet$]{$T_{m}$ est le $m$-i\`eme op\'erateur de Hecke.
Son action sur $x=(\phi\!:\!E\!\rightarrow\! E') \in  X_{0}(N)$ est
donn\'ee par $T_{m}(x)=\sum_{C}(x_{C})$, la somme portant sur tous
les sous-groupes $C$ d'ordre $m$ dans $E$ tels que $C\cap
\ker(\phi)=\{0\}$ avec $x_{C}:=(E/C\rightarrow E'/\phi(C))$.}

\vspace{0.3cm}

\item[$\bullet$]{Enfin $\sigma \in \Gal(H|k)$, avec $H$ le corps de classe de
Hilbert de $k$, correspond \textit{via} l'application d'Artin \`a la classe
d'id\'eaux $\mathcal{A}$ de $k$.}

\vspace{0.3cm}

\end{itemize}

Alors l'article de B. Gross et D. Zagier \cite{GZ} nous donne, en notant $c=c_{D}$ et $d=d_{D}$~:

\vspace{0.3cm}

\begin{center}
\begin{tabular}{lll}
$<c,T_{m}d^{\sigma}>_{\infty}$ & $\!\!\!\!=$ & $\!\!\!\!\!\displaystyle{\lim_{s\rightarrow 1}\left[-2u^{2}\sum_{n=1}^{\infty} \sigma_{\mathcal{A}}(n) r_{\mathcal{A}}(m|D|+nN) Q_{s-1}\!\left(1+\frac{2nN}{m|D|}\right)\!-\frac{h\kappa\sigma_{1}(m)}{s-1} \right]}$\\
\\
&  & $\displaystyle{+ h\kappa\left[\sigma_{1}(m)\left(\log\frac{N}{|D|}+2\sum_{p|N}\frac{\log(p)}{p^{2}-1} + 2 + 2\frac{\zeta'}{\zeta}(2)-2\frac{L'}{L}(1,\varepsilon)\right) \right]}$\\
\\
&  & $\displaystyle{+ h\kappa\left[\sum_{d|m}d\log\frac{m}{d^{2}}\right]}$\\
\\
&  & $\displaystyle{+ hur_{\mathcal{A}}(m)\left[2\frac{L'}{L}(1,\varepsilon) - 2\gamma -2\log2\pi + \log|D|\right]}$\\
\\
$<c,T_{m}d^{\sigma}>_{\mathrm{fini}}$  & $=$ & $\displaystyle{-u^{2}\sum_{1\leq n \leq m|D|/N}} \sigma_{\mathcal{A}}'(n) r_{\mathcal{A}}(m|D|-nN) \, +hur_{\mathcal{A}}(m)\log\frac{N}{m}.$\\

\end{tabular}
\end{center}

Pour les membres de droite on a~:

\vspace{0.3cm}

\begin{itemize}

\item[$\bullet$]{$r_{\mathcal{A}}(n)$ repr\'esente le nombre d'id\'eaux dans
la classe $\mathcal{A}$ de norme \'egale \`a $n$.}

\vspace{0.2cm}

\item[$\bullet$]{$\sigma_{\mathcal{A}}(n)=\sum_{d|n}\varepsilon_{\mathcal{A}}(n,d)$
et
$\sigma_{\mathcal{A}}'(n)=\sum_{d|n}\varepsilon_{\mathcal{A}}(n,d)\log\frac{n}{d^{2}}$.
On rappelle que $\varepsilon_{\mathcal{A}}(n,d)$ est nul si
$\mathrm{pgcd}(d,n/d,D)>1$. Dans le cas où $\mathrm{pgcd}(d,n/d,D)=1$, en notant $\mathrm{pgcd}(d,D)=|D_{2}|$,
$D_{1}D_{2}=D$, $\varepsilon_{D_{i}}(d)=(\frac{D_{i}}{d})$ et
$\chi_{D_{1}\!\cdot\! D_{2}}$ un certain caract\`ere de $Cl_{k}$ (voir
\cite{GZ} p. 277 et p. 268)~:
\[
\varepsilon_{\mathcal{A}}(n,d) = \varepsilon_{D_{1}}(d)\varepsilon_{D_{2}}\Big(\!-N\frac{n}{d}\Big)\chi_{D_{1}\!\cdot\! D_{2}}(\mathcal{A}).
\]
}

\item[$\bullet$]{$h=h_{k}$ est le nombre de classes associ\'e \`a $k$ et
$D=D_{k}$ est son discriminant. De plus $u=u_{k}$ est la moiti\'e du nombre de
ses unit\'es. On sait que $u=1$ sauf dans les cas $D=-3$ o\`u
$u=3$ et $D=-4$ o\`u $u=2$.}

\vspace{0.2cm}

\item[$\bullet$]{$\displaystyle{\kappa=\kappa_{N}=-12/\left(N\prod_{p|N}\Big(1+\frac{1}{p}\Big)\right)}$.}

\vspace{0.2cm}

\item[$\bullet$]{$\displaystyle{\sigma_{1}(m)=\sum_{d|m}d}$.}

\vspace{0.2cm}

\item[$\bullet$]{$\gamma\simeq 0.57$ la constante d'Euler.}

\vspace{0.2cm}

\item[$\bullet$]{$\displaystyle{\zeta(s)=\sum_{n\geq 1}\frac{1}{n^{s}}}$ est la fonction zêta de Riemann.}

\vspace{0.2cm}

\item[$\bullet$]{$\displaystyle{L(s,\varepsilon)=\sum_{n\geq 1}\frac{\varepsilon(n)}{n^{s}}}$ est la fonction $L$ de Dirichlet, avec
$\varepsilon(n)=\displaystyle{\Big(\frac{n}{D}\Big)}$.}

\vspace{0.2cm}

\item[$\bullet$]{$Q_{s-1}(t)$ est la fonction de Legendre de seconde
esp\`ece. On a plusieurs expressions de cette quantit\'e spectrale
(\cite{GZ} p. 238), par exemple pour $t>1$ et $s>0$~:
\[
Q_{s-1}(t)=\int_{0}^{\infty}\frac{\mathrm{d} u}{(t+ \sqrt{t^{2}-1}\cosh(u))^{s}}.
\]
On utilisera dans la suite la fonction
$g_{s}(z,w)=-2Q_{s-1}\left(1+\frac{|z-w|^{2}}{2\Ima(z)\Ima(w)}\right)$.
C'est en particulier une fonction holomorphe de la variable $s$ sur le domaine
$\Ree(s)>1$. Ses propri\'et\'es sont d\'etaill\'ees dans
\cite{GZ} p. 239.}
\end{itemize}

\subsubsection{Particularisations}
 La premi\`ere partie du travail consiste \`a \'evaluer cette
formule pour se ramener \`a l'expression de la hauteur du point $c
\in J_{0}(N)(H)$.

\vspace{0.3cm}

 Tout d'abord par le th\'eor\`eme de Manin-Drinfeld, $c$ et $d$
repr\'esentent la m\^eme classe dans
$J_{0}(N)(H)\otimes{\mathbb{Q}}$. On en d\'eduit l'\'egalit\'e
suivante~:
\[
<c,T_{m}d^{\sigma}>\, =\, <c,T_{m}c^{\sigma}>.
\]

On prend de plus $m=1$. On obtient alors $T_{1}c=c$. Enfin on prend $\sigma = \mathrm{Id} \in \Gal(H|k)$, ce qui impose donc de prendre $\mathcal{A}=\mathcal{O}_{k}$ l'anneau des entiers du corps $k$.
Ceci \'etant pos\'e on calcule alors le membre de droite pour
obtenir, $\widehat{h}_{J_{0}(N)}$ \'etant la hauteur de N\'eron-Tate
associ\'ee au diviseur $2\Theta$ de la vari\'et\'e ab\'elienne
$J_{0}(N)$~:

\begin{prop}\label{calcul hauteur}
\[
\widehat{h}_{J_{0}(N)}(c)=\,<c,c>\,=\,<c,c>_{\infty}\,+\,<c,c>_{\mathrm{fini}},
\]

avec~:

\begin{tabular}{llll}
$<c,c>_{\infty}$ & $\!\!\!=$ & $\!\!\!\!\displaystyle{\lim_{s\rightarrow 1}\left[-2u^{2}\sum_{n=1}^{\infty} \sigma_{\mathcal{O}_{k}}(n) r_{\mathcal{O}_{k}}(|D|+nN) Q_{s-1}\!\left(1+\frac{2nN}{|D|}\right)\!-h\frac{\kappa}{s-1} \right]}$ & $\!\!(i)$\\
\\
&  & $\displaystyle{+ h\kappa\left(\log\frac{N}{|D|}+2\sum_{p|N}\frac{\log(p)}{p^{2}-1} + 2 + 2\frac{\zeta'}{\zeta}(2)-2\frac{L'}{L}(1,\varepsilon)\right)}$ & $\!\!(ii)$\\
\\
&  & $\displaystyle{+ hu\left[2\frac{L'}{L}(1,\varepsilon) - 2\gamma -2\log2\pi + \log|D|\right]}$ & $\!\!\!(iii)$\\
\\
$<c,c>_{\mathrm{fini}}$  & $=$ & $\displaystyle{-u^{2}\sum_{1\leq n \leq |D|/N}} \sigma_{\mathcal{O}_{k}}'(n) r_{\mathcal{O}_{k}}(|D|-nN) \, +hu\log(N)$ & $\!\!(iv)$\\

\end{tabular}

\end{prop}

\vspace{0.3cm}

\subsection{Preuve du théorème \ref{equiv}}
On se place toujours dans le m\^eme cadre, le discriminant $D$ du
corps $k$ est fix\'e avec les conditions de l'introduction.
D'apr\`es la proposition \ref{NK} l'ensemble $\mathbb{N}_{k}$ est
infini, on peut donc faire tendre $N$ vers l'infini. On s'efforce
alors dans cette troisi\`eme partie de trouver un \'equivalent,
lorsque $N$ tend vers l'infini, de la hauteur
$\widehat{h}_{J_{0}(N)}(c)$. Nous allons donc \'etudier la
contribution de chaque terme de la proposition \ref{calcul
hauteur}. On commence par donner quelques majorations utiles.

\subsubsection{Majorations}

\begin{lem}\label{lem1n}
Si on note $\tau(n)$ le nombre de diviseurs de $n$, alors on a les
majorations $|\sigma_{\mathcal{O}_{k}}(n)|\leq \tau(n)$ et
$|\sigma_{\mathcal{O}_{k}}'(n)|\leq \tau(n)\log(n)$.
\end{lem}

\begin{proof} Il suffit de voir que
$|\varepsilon_{\mathcal{A}}(n,d)|$ est born\'e par $1$.
\end{proof}
\vspace{0.3cm}

\begin{lem}\label{lem2n}
On rappelle la majoration~:
$\tau(n)=O_{\varepsilon}(n^{\varepsilon})$ pour tout $\varepsilon
>0$.
\end{lem}

\begin{proof} On se reportera à \cite{Ten} et \cite{Ten2} p. 13 et suivantes.
\end{proof}
\vspace{0.3cm}

\begin{lem}\label{lem3n}
On peut majorer~: $r_{\textit{O}_{k}}(n)
=O_{\varepsilon}(n^{\varepsilon})$ pour tout $\varepsilon
>0$.
\end{lem}

\begin{proof} Dans un anneau d'entiers, un id\'eal se
d\'ecompose en produit d'id\'eaux premiers. Soient $n\geq 1$ et
$\mathcal{I}=\mathcal{P}_{1}^{\alpha_{1}}...\mathcal{P}_{l}^{\alpha_{l}}$
un id\'eal de norme $n$. En prenant la norme on obtient une
\'egalit\'e du type $n=p_{1}^{\beta_{1}}...p_{l}^{\beta_{l}}$
avec les $p_{i}$ des entiers naturels premiers. Les possibilit\'es
pour l'id\'eal $\mathcal{I}$ sont donc fonction du nombre
d'id\'eaux au-dessus de chaque premier $p_{i}|n$. Puisque $k$ est
un corps quadratique, il y a au plus deux id\'eaux au-dessus d'un
entier premier $p$ de $\mathbb{Z}$ (auquel cas $p$ est totalement
d\'ecompos\'e), ceci donne donc lieu \`a au plus $2^{l}$ id\'eaux
$\mathcal{I}$ de norme $n$. Or $l=\sum_{p|n}1=\omega(n)$ et par d\'efinition de $\tau(n)$ on a~: $2^{\omega(n)}\leq \tau(n)$. On
conclut donc par le lemme pr\'ec\'edent.
\end{proof}
\vspace{0.3cm}

\begin{lem}\label{lem4n}
Si $s>1$ on a les propri\'et\'es asymptotiques suivantes~:
\begin{align*} 
Q_{s-1}(t)=&\, O_{t\rightarrow +\infty}(t^{-s}),\\
Q_{s-1}(t)=&-\frac{1}{2}\log(t-1)+O_{t\rightarrow 1}(1).
\end{align*}
\end{lem}

\begin{proof} On pourra par exemple se r\'ef\'erer \`a
\cite{EMOT} à partir de la page 155.
\end{proof}

\begin{lem}\label{lem6n}
On rappelle enfin~: $\displaystyle{\zeta(s)=\frac{1}{s-1}+O_{s\rightarrow
1}(1)}$.
\end{lem}

\begin{proof} Il suffit de faire une comparaison
s\'erie-int\'egrale.
\end{proof}

\subsubsection{Les termes (ii), (iii), (iv)}
 Le traitement des trois derniers termes est assez rapide. En effet
\`a $D$ fix\'e $h=h_{k}$ est constant, on obtient donc directement
que le terme (iii) est un $O_{N\rightarrow +\infty}(1)$. De plus
en utilisant l'estimation aisée $\kappa=O\left(\frac{1}{N}\right)$ on a imm\'ediatement que (ii) est
un $O\left(\frac{\log(N)}{N}\right)$.

\vspace{0.3cm}

 Pour (iv) on remarque que $\displaystyle{-u^{2}\sum_{1\leq n \leq |D|/N}}
\sigma_{\mathcal{O}_{k}}'(n) r_{\mathcal{O}_{k}}(|D|-nN) = 0$ si
$N>|D|$ ; cela suffit puisqu'on va considérer $N$ grand à $D$ fixé. On obtient donc que le terme dominant est $hu\log(N)$.

\vspace{0.3cm}

Jusqu'ici on a donc montr\'e que le terme principal des
contributions (ii), (iii) et (iv) est $h\,u\,\log(N)$. Il nous
reste maintenant \`a \'etudier le terme (i) issu (comme (ii) et
(iii)) des places archim\'ediennes.

\subsubsection{Le terme (i)}

Le but de toute cette partie est de g\'en\'eraliser la d\'emarche
de P. Michel et E. Ullmo dans \cite{MUll} pour montrer une
majoration du terme (i) par un $O_{\varepsilon}(N^{\varepsilon-1})$.
Commen\c{c}ons par poser~:
\[
\widetilde{H}(s):=\displaystyle{-2u^{2}\sum_{n=1}^{\infty} \sigma_{\mathcal{O}_{k}}(n)
r_{\mathcal{O}_{k}}(|D|+nN) Q_{s-1}\left(1+\frac{2nN}{|D|}\right).}
\]

Les lemmes \ref{lem2n}, \ref{lem3n} et \ref{lem4n} assurent que $\widetilde{H}$ converge absolument pour $\Ree(s)>1$ et d\'efinit une fonction holomorphe.

On va utiliser plusieurs fonctions introduites dans l'article de
B. Gross et D. Zagier pour \'etudier la fonction $\widetilde{H}$ au
voisinage de $1$. La fonction $g_{s}$ est d\'efinie en $2.1$. On
rappelle donc~:

Pour $z,z'\in\mathcal{H}$, on pose (\cite{GZ} p. 251 et 252)~:
\[
G_{N,s}^{1}(z,z'):=\sum_{\substack{\gamma\in
\Gamma_{0}(N)/\{\pm 1\} \\ \gamma z'\neq z}}g_{s}(z,\gamma
z') + 2hu\left[\frac{\Gamma'}{\Gamma}(s)
-\log(2\pi)+\frac{L'}{L}(1,\varepsilon)+\frac{1}{2}\log(|D|)\right].
\]

Pour $\mathcal{A}\in Cl_{k}$, on pose (\cite{GZ} p. 243)~:
\[
\gamma_{N,s}^{1}(\mathcal{A}):=\sum_{\substack{\mathcal{A}_{1},\mathcal{A}_{2}\in Cl_{k} \\\mathcal{A}_{1}\mathcal{A}_{2}^{-1}=\mathcal{A}}}G_{N,s}^{1}(\tau_{\mathcal{A}_{1}}, \tau_{\mathcal{A}_{2}}).
\]

Les calculs de B. Gross et D. Zagier montrent alors que (\cite{GZ}
p. 243 et p. 247 combin\'ee avec p. 285)~:
\[
\widetilde{H}(s)= \gamma_{N,s}^{1}(\mathcal{O}_{k})=\sum_{\substack{\mathcal{A}_{1},\mathcal{A}_{2}\in Cl_{k} \\\mathcal{A}_{1}\mathcal{A}_{2}^{-1}=\mathcal{O}_{k}}}G_{N,s}^{1}(\tau_{\mathcal{A}_{1}}, \tau_{\mathcal{A}_{2}})= \sum_{\mathcal{A}_{1}\in Cl_{k}}G_{N,s}^{1}(\tau_{\mathcal{A}_{1}}, \tau_{\mathcal{A}_{1}}).
\]

On peut alors trouver dans l'article de P. Michel et
E. Ullmo (\cite{MUll} p. 673) une \'etude d'un terme
$G_{N,s}^{1}(\tau, \tau)$ ($=H(s)$ dans leur notation) pour un
point de Heegner $\tau$. Leur r\'esultat est le suivant~:

\begin{prop}\label{Phrag-Lind}
On a la majoration suivante, valable pour tout $\varepsilon >0$~:
\[
\lim_{s\rightarrow 1}\left(G_{N,s}^{1}(\tau,
\tau)-\frac{\kappa}{s-1}\right)=O_{\varepsilon}(N^{\varepsilon-1}).
\]
\end{prop}

\begin{proof} 

Nous allons suivre \cite{MUll} en remarquant que
leur preuve reste valide pour tout point de Heegner $\tau$ lorsque
$D$ est fix\'e. On introduit le noyau automorphe pour
$\Gamma_{0}(N)/\{\pm 1\}$~:
\[
\displaystyle{G_{N,s}(z,z')=\sum_{\substack{\gamma\in \Gamma_{0}(N)/\{\pm
1\}\\ \gamma z'\neq z}}g_{s}(z,\gamma z')}.
\]
On a alors l'\'egalit\'e pour $a>1$~:

\begin{equation}\label{GnsmoinsGna}
G_{N,s}^{1}(\tau, \tau)-G_{N,a}^{1}(\tau, \tau)=G_{N,s}(\tau,\tau)-G_{N,a}(\tau,\tau)-2hu\left(\frac{\Gamma'}{\Gamma}(s)-\frac{\Gamma'}{\Gamma}(a)\right).
\end{equation}

On utilise pour conclure le lemme la proposition suivante, dont la preuve figure dans le livre de H. Iwaniec \cite{Iw} p. 105, th\'eor\`eme $7.5$. 

\begin{prop}
Soient $a>1$ et $\Ree(s)>1$. Soit $N$ un entier sans facteur carré et premier à $6$. On note $s_{j}(1-s_{j})$ la $j$-i\`eme
valeur propre du laplacien sur $X_{0}(N)$. On prend de plus
$(u_{j})_{j}$ une base orthonormale de fonctions propres
associ\'ees \`a ces valeurs propres. On note de plus $E_{\rho}$ la
s\'erie d'Eisenstein associ\'ee \`a la pointe $\rho$. On pose enfin~:
\[
\chi_{sa}(v)=\frac{1}{(s-v)(1-s-v)}-\frac{1}{(a-v)(1-a-v)}.
\]

On a alors l'\'egalit\'e~:

\vspace{0.1cm}

\begin{tabular}{lll}
$\displaystyle{G_{N,s}(z,z')-G_{N,a}(z,z')}$ & $=$ & $\displaystyle{\sum_{j}\chi_{sa}(s_{j})u_{j}(z)\overline{u}_{j}(z')}$\\

&  & $\displaystyle{+\sum_{\rho\in\{\mathrm{Pointes}\}}\frac{1}{4\pi i}\int_{\frac{1}{2}+i\mathbb{R}}\chi_{sa}(v)E_{\rho}(z,v)\overline{E}_{\rho}(z',\overline{v})dv ,}$\\

\end{tabular}

et la s\'erie et l'int\'egrale convergent absolument et
uniform\'ement sur tout compact.
\end{prop}

\end{proof}

Cette proposition permet de prolonger la fonction
$G_{N,s}(z,z')-G_{N,a}(z,z')$ en une fonction m\'eromorphe
(que l'on notera de la m\^eme mani\`ere) sur le domaine
$\Re(s)>1/2$ avec un p\^ole simple en $s=1$ de r\'esidu \'egal \`a
$\kappa$ (voir aussi \cite{GZ} p. 239). Cette fonction est m\^eme
holomorphe sur le domaine $\Ree(s)>3/4$ priv\'e du point $s=1$
(voir \cite{MUll} p. 672).

On va choisir une bande verticale dans le plan complexe contenant
l'abscisse $s=1$ dans le but d'appliquer le principe de
Phragmen-Lindelöf \`a la fonction~:
\[
G_{N,s}(\tau,\tau)-G_{N,a}(\tau,\tau) -\frac{\kappa}{s-1}.
\]
On en d\'eduira une majoration au voisinage de $s=1$ de celle-ci et donc
par l'\'egalit\'e (\ref{GnsmoinsGna}) de $G_{N,s}^{1}(\tau, \tau)-\kappa/(s-1)$.

\vspace{0.3cm}

Tout d'abord par application des lemmes \ref{lem2n}, \ref{lem3n},
\ref{lem4n} et \ref{lem6n} on a la majoration valable pour
$\Re(s)>1$ et $\varepsilon>0$~:
\[
G_{N,s}^{1}(\tau, \tau)=O_{\varepsilon}\left(N^{\varepsilon
-1}\left(1+\frac{1}{|s-1|}\right)\right).
\]

On va \`a pr\'esent \'etudier la croissance de la fonction
$G_{N,s}^{1}(\tau, \tau)-\kappa/(s-1)$ (toujours par
l'interm\'ediaire de l'\'egalit\'e (\ref{GnsmoinsGna}) ci-avant) dans une bande du type
$b\leq s\leq 1+\varepsilon$ avec $3/4 < b < 1$. On prend par
exemple $b=7/8$. On utilise alors $t=\Ima(s)$ et on a~:

\begin{prop}
Il existe un $A>0$ tel que pour tout $s$ dans la bande $7/8\leq
s\leq 1+\varepsilon$ on a la majoration, avec $t=\Ima(s)$ et $\tau$
un point de Heegner~:
\[
G_{N,s}(\tau,\tau)-G_{N,a}(\tau,\tau)-\frac{\kappa}{s-1}\ll(N(1+|t|))^{A},
\]
et la constante implicite est ind\'ependante de $N$.
\end{prop}

\begin{proof} On se reportera \`a \cite{MUll} p. 673.
\end{proof}

On a donc une croissance suffisamment faible de la fonction
$G_{N,s}^{1}(\tau, \tau) -\kappa/(s-1)$ dans la bande $7/8\leq
\Ree(s) \leq 1+\varepsilon$ pour pouvoir conclure par le principe
de Phragmen-Lindelöf en utilisant l'\'equation fonctionnelle de
\cite{MUll} p. 654. Ceci ach\`eve la preuve du th\'eor\`eme
\ref{Phrag-Lind}.

Il suffit alors d'appliquer ce th\'eor\`eme \`a chacune des $h$
classes associ\'ees \`a $k$ pour obtenir~:

\begin{prop}\label{terme difficile}
On a la majoration valable pour tout $\varepsilon>0$~:
\[
\lim_{s\rightarrow 1}\left(\widetilde{H}(s)-h\frac{\kappa}{s-1}\right)=O_{\varepsilon}(N^{\varepsilon-1}).
\]
\end{prop}

La conjonction de la proposition \ref{calcul hauteur} avec les résultats de majorations \ref{lem1n} à \ref{lem6n} et \ref{terme difficile} démontre le théorème \ref{equiv}.

\subsection{Corollaires}

On donne ici deux corollaires du théorème \ref{equiv}~:

\begin{cor}
Sous les hypoth\`eses du th\'eor\`eme, avec $|D|>4$, on a
l'estimation de $h_{k}$ suivante~:
\[
\lim_{\substack{N\rightarrow \infty \\N\in \mathbb{N}_{k}}}\frac{\widehat{h}_{J_{0}(N)}(c_{D})}{\log(N)} =h_{k}.
\]

\end{cor}

\begin{cor}\label{grandN}
Pour tout discriminant $D\equiv 1 \pmod{4}$ n\'egatif et sans
facteur carr\'e il existe un entier $N_{0}(D)$ et un certain
ensemble $\mathbb{N}_{k}$ de congruences modulo $D$ tels que pour
tout $N\geq N_{0}(D)$ et $N\in \mathbb{N}_{k}$ les points de
Heegner $c_{D}$ sont d'ordre infini dans $J_{0}(N)(H)$.
\end{cor}

\begin{remarque}
Pour certains petits $N$ on sait que le corps de d\'efinition est n\'ecessairement strictement plus gros que $\mathbb{Q}$. Par exemple on a
$c_{D}\in{J_{0}(N)(H)\backslash J_{0}(N)(\mathbb{Q})}$ pour tout
$N\in{\{11,14,15,17,19,20,21,24,27,32,36,49\}}$. En effet on
conna\^it les courbes modulaires $X_{0}(N)$ de genre 1, il y en a
douze et elles correspondent aux niveaux $N$ ci-avant. Ces courbes sont
des courbes elliptiques et on a donc dans ce cas un isomorphisme entre
$J_{0}(N)$ et $X_{0}(N)$. On sait de plus que pour ces douze courbes le
groupe $X_{0}(N)(\mathbb{Q})$ est fini, i.e. leur rang sur
$\mathbb{Q}$ est \'egal \`a z\'ero. (On se base ici sur une
\'etude des courbes modulaires de genre 1 dont la r\'ef\'erence
est \cite{Lig}.)
\end{remarque}

\subsection{Remarques sur la conjecture de Lang et Silverman}

On rappelle qu'on note $\hF(A/k)$ pour la hauteur de Faltings (relative)
d'une vari\'et\'e ab\'elienne $A$. La hauteur de Néron-Tate associée à $2\Theta$ sera notée $\widehat{h}_{A}(.)$.
Rappelons l'énoncé de la conjecture de Lang et Silverman sur les variétés abéliennes~:

 \begin{conj}(Lang, Silverman)
Soit $g\geq 1$ et soit $k$ un corps. Il existe une constante
$c=c(k,g)>0$ telle que pour toute vari\'et\'e ab\'elienne $A/k$ de
dimension $g$, tout diviseur ample et symétrique $\mathcal{D}\in \Divi(A)$ et tout point
$P \in A(k)$ tel que $\overline{{\mathbb{Z}\!\cdot\! P}}=A$ on a~:
\[
\widehat{h}_{A,D}(P)\geq c\; \hF(A/k).
\]
\end{conj}

On trouve dans l'article \cite{JoKr} l'\'equivalent suivant, obtenu par des techniques de géométrie hyperbolique complexe~:

\begin{thm}(Jorgenson, Kramer)
Pour tout $N$ sans facteur carr\'e et premier à $6$ on a, lorsque $N$ tend vers l'infini~:
\[
\hst(J_{0}(N))=\frac{g(N)}{3}\log(N)\,+\,o\Big(g(N)\log(N)\Big).
\]
\end{thm}

Or la dimension de $J_{0}(N)$ est \'egale au genre de $X_{0}(N)$
et ce dernier s'exprime comme suit lorsque $N$ est sans facteur
carr\'e (voir \cite{Shi})~:
\[
g(X_{0}(N))=1+\frac{N}{12}\prod_{p|N}\left(1+\frac{1}{p}\right)-\frac{1}{4}\prod_{p|N}\left(1+\left(\frac{-1}{p}\right)\right)-\frac{1}{3}\prod_{p|N}\left(1+\left(\frac{-3}{p}\right)\right)-\frac{1}{2}\tau(N),
\]

et donc en particulier on a l'\'equivalent lorsque $N\rightarrow
+\infty$~:

\[
g(X_{0}(N))\sim \frac{N}{12}\prod_{p|N}\left(1+\frac{1}{p}\right),
\]

avec l'encadrement pour une certaine constante $\alpha>0$~:
\[
1\leq \prod_{p|N}\left(1+\frac{1}{p}\right) \leq \alpha\log\log(N).
\]

\begin{remarque}(importante)
Les points de Heegner v\'erifient la
condition $\overline{{\mathbb{Z}\!\cdot\! P}}=J_{0}(N)$ au moins lorsque le
discriminant $D$ est choisi suffisamment grand. Cela d\'ecoule de
l'\'etude men\'ee par J. Nekov\'ar et N. Schappacher dans
\cite{Nek}.
\end{remarque}

On a donc~:

\begin{fait} 
Si nous r\'eunissons ici les r\'esultats
obtenus dans les th\'eor\`emes 3.11 et 3.17~:
$\widehat{h}_{J_{0}(N)}(c)\sim h_{k}u_{k}\log(N)$ et $\hst(J_{0}(N))\sim
g(N)\log(N)/3$, on voit que la conjonction de ces th\'eor\`emes
constitue un exemple indiquant que la constante pr\'esente dans la
conjecture de Lang et Silverman doit n\'ecessairement d\'ependre de
la dimension $g$ des vari\'et\'es ab\'eliennes consid\'er\'ees. On
peut m\^eme affirmer pour un corps quadratique imaginaire $k$
donn\'e et $H$ son corps de classe de Hilbert~:
\[
c(H,g)\leq \frac{3h_{k}}{g}=\frac{3[H:k]}{g}.
\]

\end{fait}

On déduit de plus de la comparaison des deux asymptotiques le corollaire \ref{asympt}.

\subsection{Ordre et niveau sur la jacobienne}

On compare ici le r\'esultat asymptotique montr\'e
pr\'ec\'edemment avec un r\'esultat existant pour les petites
valeurs du niveau $N$. On pourra se r\'ef\'erer \`a l'article de
H. Nakazato \cite{Na}. On se place dor\'enavant dans la situation
suivante~: soit $E$ une courbe elliptique d\'efinie sur
$\mathbb{Q}$. D'apr\`es le th\'eor\`eme de Wiles \'etendu par
Breuil, Conrad, Diamond, Taylor elle est munie d'un
morphisme non constant~:
\[
\varphi:X_{0}(N)\longrightarrow E(\mathbb{C}).
\]

On note de plus $E_{l}$ l'ensemble des points de $l$-torsion de $E$. Si $E$ n'est pas \`a multiplication complexe on pose~:
\[
S_{E}=\Big\{l\,|\,\Gal(\mathbb{Q}(E_{l})|\mathbb{Q})\neq \Aut(E_{l})\Big\}\cup \{l|N\}\cup \{2,3\}.
\]

Si $E$ est \`a multiplication complexe on posera seulement~:
\[
S_{E}=\{l|N\}\cup \{2,3\}.
\]

Fixons alors un discriminant $D\equiv 1 \pmod{4}$ n\'egatif tel qu'aucun facteur premier de $D$ ne soit inclu dans $S_{E}$.
L'ensemble $S_{E}$ est fini et il y a une
infinit\'e de tels $D$ (voir \cite{Na}). Soit alors $N\in
\mathbb{N}_{k}$. Dans ces conditions on a~:

\begin{thm}(Nakazato)\label{petitN}
Soit $\tau \in X_{0}(N)$ un point de Heegner associ\'e \`a $k$. Si on a
$h_{k}>\deg(\varphi)$ alors $\varphi(\tau)$ est un point d'ordre
infini sur $E$.
\end{thm}

\begin{remarque}
 On peut minorer le degr\'e de $\varphi$ en
fonction de $N$, par exemple en suivant \cite{Wat}~:
\[
\deg(\varphi)\geq N^{\frac{7}{6}-\varepsilon}.
\]
\end{remarque}

En utilisant cette remarque et les th\'eor\`emes \ref{grandN} et
\ref{petitN} on pourra garder à l'esprit que les points de Heegner sur la jacobienne $J_{0}(N)$ sont génériquement des points d'ordre infini.

\section{Variation du corps et restriction des scalaires à la Weil}

On présente dans cette partie l'effet de la variation du corps, vue sous deux angles différents, sur l'énoncé de la conjecture de Lang et Silverman.

\subsection{Variation du corps}

Nous commençons notre étude par la remarque suivante~: soient $k$ un corps de nombres et $(A,\Theta)/k$ une variété abélienne polarisée par un diviseur $\Theta$. Soit $P\in{A(k)}$ un point d'ordre infini. Pour tout $N\geq1$ posons $P_{N}=\frac{1}{N}P$ et $k_{N}=k[P_{N}]$. Alors~:
\[
\widehat{h}_{A,\Theta}(P_{N})=\frac{1}{N^{2}}\widehat{h}_{A,\Theta}(P).
\]

On a donc~:

\begin{fait}(classique)
Soient $k$ un corps de nombres et $A/k$ une variété abélienne. Il existe une suite d'extensions $k_{N}/k$ et une suite de points $P_{N}\in{A(k_{N})}$ vérifiant $\mathbb{Z}.P_{N}$ Zariski-dense telles que~:

\[
\lim_{N\rightarrow +\infty}\widehat{h}_{A,\Theta}(P_{N})=0.
\]

\end{fait}

Ce fait souligne la nécessité d'avoir une constante dépendant du corps $k$ dans l'énoncé de la conjecture de Lang et Silverman. On ajoute que la dépendance minimale en le corps $k$ devrait être en $[k:\mathbb{Q}]^{-1/g}$.

\subsection{Restriction des scalaires à la Weil}\label{restriction de weil}
Nous exploitons ici la même idée de division de points, mais présentée différemment.
Nous allons nous intéresser à la restriction des scalaires, appelée aussi foncteur norme. On consultera \cite{Mil} et \cite{Weil2} pour une définition plus générale et les preuves des propriétés utilisées ci-après.

\vspace{0.3cm}

Soient $k$ un corps de nombres et $L/k$ une extension finie de degré
$m$. Soit $A/L$ une variété abélienne définie sur $L$ de dimension
$g$. On s'intéresse à la variété obtenue par restriction des scalaires $A_{*}=\NL A$ définie sur $k$ par restriction des scalaires. C'est une variété abélienne car pour toute extension galoisienne $k'$ de $k$ contenant $L$ on a l'isomorphisme~:
\[
\psi:(\NL A)_{k'}\longrightarrow A^{\sigma_{1}}_{k'}\!\times\!...\!\times\! A^{\sigma_{m}}_{k'},
\]

où les $\sigma_{i}$ sont les plongements de $L$ dans $k'$ au-dessus de $k$. Elle arrive de plus équipée d'un morphisme surjectif~:
\[
p:(\NL A)_{L} \longrightarrow A.
\]

Soit $b\in{\Pico(A)}$. On lui associe l'élément $b_{*}=p^{\sigma_{1}*}(b^{\sigma_{1}})+...+p^{\sigma_{m}*}(b^{\sigma_{m}})$. Alors la proposition 4 de \cite{Mil} nous donne un isomorphisme~:
\begin{align*} 
\Pico(A) &\to \Pico(A_{*})\\
b & \mapsto b_{*}.
\end{align*}

On associera à la variété abélienne polarisée $(A,\Theta)/L$ sa restriction $(A_{*},\Theta_{*})/k$ \textit{via} cet isomorphisme. La proposition $5$ de \cite{Mil} nous donne alors~:

\begin{prop}(Milne)

On note $<.,.>_{L}:\Pico(A)\!\times\! A(L)\rightarrow \mathbb{R}$ l'accouplement de Néron-Tate. Soient $a\in{A_{*}(k)}$ et $b\in{\Pico(A)}$, alors~:
\[
<b_{*},a>_{k}=<b,p(a)>_{L}.
\]

\end{prop}

Considérons l'application (surjective, voir par exemple \cite{HiSi} p. 208) suivante~:
\begin{align*} 
\Phi_{\Theta} \colon A &\to \Pico(A)\\
Q & \mapsto t_{Q}^{*}\Theta-\Theta.
\end{align*}

Nous choisissons alors dans l'énoncé de Milne~: $b=\Phi_{\Theta}(p(a))$, image du point $p(a)\in{A}$ dans $\Pico(A)$. On sait alors que $<\Phi_{\Theta}(p(a)),p(a)>_{L}=\widehat{h}_{A/L,\Theta}(p(a))$ et on déduit de ce qui précède~:
\[
\widehat{h}_{A_{*}/k,\Theta_{*}}(a)=\widehat{h}_{A/L,\Theta}(p(a)).
\]

Considérons alors le cas de figure suivant. On choisit un point $P_{1}\in{A(k)}$ d'ordre infini et on forme, pour $N\geq 1$, une suite de points $P_{N}=\frac{1}{N}P_{1}\in{A(k_{N})}$, avec $k_{N}=k[P_{N}]$. On note $m_{N}=[k_{N}:k]$. On peut donc définir une suite de variétés abéliennes $A_{N}=\NkN A$ définies sur $k$ telles que, en notant $P_{N}'$ un antécédent de $P_{N}$ par $p$ dans $A_{N}(k)$~:
\[
\widehat{h}_{A_{N}/k,\Theta_{N}}(P_{N}')=\widehat{h}_{A/k_{N},\Theta}(P_{N})=\frac{1}{N^{2}}\widehat{h}_{A/k_{N},\Theta}(P_{1})=\frac{1}{N^{2}}\widehat{h}_{A/k,\Theta}(P_{1}).
\]

Dans le même temps, grâce à l'isomorphisme $\psi$ on peut déduire les relations suivantes sur les dimensions et les hauteurs stables~:
\[
\dime(A_{N})=m_{N}\dime(A),\;\;\;\;\; h_{st}(A_{N})=m_{N}h_{st}(A).
\]

On doit enfin calculer l'adhérence de $\mathbb{Z}\!\cdot\! P_{N}'$. Or comme $P_{N}=\frac{1}{N}P_{1}$, on aura  $\mathbb{Z}\!\cdot\! P_{N}=\frac{1}{N}\mathbb{Z}\!\cdot\! P_{1}$, donc~:
\[
\overline{\mathbb{Z}\!\cdot\! P_{N}'}= \Big\{(P,...,P)\in{A^{\sigma_{1}}\!\times\!...\!\times\! A^{\sigma_{m_{N}}}}\Big\}\subsetneq A^{\sigma_{1}}\!\times\!...\!\times\! A^{\sigma_{m_{N}}}
\]

On peut donc garder à l'esprit~:

\begin{fait}

Soit $k$ un corps de nombres. Il existe une suite de variétés abéliennes $A_{N}/k$ définies sur $k$ et une suite de points d'ordre infini $Q_{N}\in{A_{N}(k)}$ telles que $(\dime(A_{N}))_{N\geq 1}$ est croissante et~:
\begin{align*} 
&\lim_{N\rightarrow +\infty}\widehat{h}_{A_{N}}(Q_{N})= 0,\\
&\lim_{N\rightarrow +\infty}h_{st}(A_{N})= +\infty,\\
&\dime(\overline{\mathbb{Z}.Q_{N}})= \dime(A_{1}).
\end{align*}
\end{fait}

Ce fait souligne le caractère crucial de l'hypothèse $\overline{\mathbb{Z}\!\cdot\! P}=A$ dans l'énoncé de la conjecture de Lang et Silverman. 

\begin{remarque}
Une variante consiste à considérer la situation $A=A_{1}\!\times\! A_{2}$ et un point $P=(P_{1},O)\in{A(k)}$, avec $\hF(A_{2}/k)$ très grand.
\end{remarque}

\end{document}